\documentclass[english]{amsart}

\usepackage{amsmath}

\usepackage{amssymb}

\usepackage{amscd} \usepackage{tabularx}
\usepackage[all,cmtip,line]{xy}

\newcommand{\ra}{\rightarrow}
\newcommand{\lra}{\longrightarrow}
\newcommand{\la}{\leftarrow}
\newcommand{\lla}{\longleftarrow}
\newcommand{\into}{\hookrightarrow}

\newcommand{\da}{\downarrow}
\newcommand{\iso}{\stackrel{\sim}{\ra}}
\newcommand{\liso}{\stackrel{\sim}{\lra}}

\newcommand{\pfbegin}{{{\em Proof:}\;}}
\newcommand{\pfend}{$\Box$ \medskip}

\newlength{\ownl}

\newcommand{\ndiv}{{\mbox{$\not| $}}}

\newcommand{\Aff}{{\operatorname{Aff}}}

\newcommand{\Aut}{{\operatorname{Aut}\,}}

\newcommand{\Fil}{{\operatorname{Fil}\,}}

\newcommand{\Frob}{{\operatorname{Frob}}}

\newcommand{\Gal}{{\operatorname{Gal}\,}}

\newcommand{\Hom}{{\operatorname{Hom}\,}}

\renewcommand{\Im}{{\operatorname{Im}\,}}

\newcommand{\DR}{{\operatorname{DR}}}
\newcommand{\WD}{{\operatorname{WD}}}

\newcommand{\gr}{{\operatorname{gr}\,}}

\newcommand{\rec}{{\operatorname{rec}}}
\newcommand{\tr}{{\operatorname{tr}\,}}

\newcommand{\ann}{{\operatorname{an}}}

\newcommand{\Fsemis}{{\operatorname{F-ss}}}

\newcommand{\A}{{\mathbb{A}}}

\newcommand{\C}{{\mathbb{C}}}

\newcommand{\G}{{\mathbb{G}}}

\newcommand{\PP}{{\mathbb{P}}}
\newcommand{\Q}{{\mathbb{Q}}}
\newcommand{\R}{{\mathbb{R}}}

\newcommand{\Z}{{\mathbb{Z}}}

\newcommand{\CH}{{\mathcal{H}}}

\newcommand{\CL}{{\mathcal{L}}}
\newcommand{\CM}{{\mathcal{M}}}

\newcommand{\CR}{{\mathcal{R}}}

\newcommand{\barF}{\overline{{F}}}

\newcommand{\barK}{\overline{{K}}}

\newcommand{\barM}{\overline{{M}}}

\newcommand{\barQQ}{\overline{{\Q}}}

\newcommand{\barr}{\overline{{r}}}

\newcommand{\tc}{\widetilde{{c}}}

 \newcommand{\tsigma   }{\widetilde{\sigma}}   
    
 \newcommand{\ttau     }{\widetilde{\tau}}

\def\RCS$#1: #2 ${\expandafter\def\csname RCS#1\endcsname{#2}}
\RCS$Revision: 364 $
\RCS$Date: 2012-09-24 15:39:31 -0400 (Mon, 24 Sep 2012) $

\newcommand{\cL}{\mathcal{L}}

\newcommand{\Favoid}{F^{({avoid})}}

\newcommand{\mc}{\mathcal}

\DeclareMathOperator{\GO}{GO}

\newcommand{\GL}{\operatorname{GL}}
\newcommand{\GSp}{\operatorname{GSp}}

\newcommand{\HT}{\operatorname{HT}}

\newcommand{\cF}{\mathcal{F}}

\newcommand{\cH}{\mathcal{H}}

\newcommand{\cM}{\mathcal{M}}

\newcommand{\cO}{\mathcal{O}}

\newcommand{\barFil}{{\overline{\Fil}}}
\newcommand{\bargr}{{\overline{\gr}}}

\usepackage{hyperref}

\usepackage{amsthm}

 \newtheorem{ithm}{Theorem}
\newtheorem{icor}[ithm]{Corollary}
\newtheorem{thm}{Theorem}[section]

\newtheorem{cor}[thm]{Corollary}
 
\newtheorem{lem}[thm]{Lemma} 
 \theoremstyle{definition}
 \theoremstyle{definition}
 \theoremstyle{remark}

\numberwithin{equation}{section}

\theoremstyle{definition}

\begin{document}
\title{Automorphy and irreducibility of some $l$-adic representations.}

\author{Stefan Patrikis}\email{patrikis@ias.edu}\address{Institute for Advanced Study, Princeton} \thanks{The  first author was partially
  supported by NSF grant DMS-1062759}
\author{Richard Taylor}\email{rtaylor@ias.edu}\address{Institute for Advanced Study, Princeton} \thanks{The  second author was partially
  supported by NSF grant DMS-1252158}  
\begin{abstract} In this paper we prove that a pure, regular, totally odd, polarizable weakly compatible system of $l$-adic representations is potentially automorphic. The innovation is that we make no irreducibility assumption, but we make a purity assumption instead. For compatible systems coming from geometry, purity is often easier to check than irreducibility. We use Katz's theory of rigid local systems to construct many examples of motives to which our theorem applies. We also show that if $F$ is a CM or totally real field and if $\pi$ is a polarizable, regular algebraic, cuspidal automorphic representation of $GL_n(\A_F)$, then for a positive Dirichlet density set of rational primes $l$, the $l$-adic representations $r_{l,\imath}(\pi)$ associated to $\pi$ are irreducible.
\end{abstract}
\maketitle
\newpage

\section*{Introduction.}\label{sec:intro}

This note is a postscript to \cite{blggt}. 

Suppose that $F$ and $M$ are
number fields, that $S$ is a finite set of primes of $F$ and that $n$
is a positive integer. By a {\em weakly compatible system} of $n$-dimensional $l$-adic representations of $G_F$ defined over $M$ and unramified outside $S$ we shall mean a family
of continuous semi-simple representations
\[ r_\lambda: G_F \lra \GL_n(\barM_\lambda), \]
where $\lambda$ runs over the finite places of $M$, with the following properties.
\begin{itemize}
\item If $v\notin S$ is a finite place of $F$, then for all $\lambda$ not dividing the residue characteristic of $v$, the representation $r_\lambda$ is unramified at $v$ and the characteristic
polynomial of $r_\lambda(\Frob_v)$ lies in $M[X]$ and is independent of $\lambda$. 
\item Each
representation $r_\lambda$ is de Rham at all places above the residue characteristic of $\lambda$, and in fact crystalline at any place $v \not\in S$ which divides
the residue characteristic of $\lambda$. 
\item For each embedding $\tau:F \into \barM$ the $\tau$-Hodge--Tate numbers of $r_\lambda$ are independent of $\lambda$.
\end{itemize}

In this paper we prove the following theorem (see Theorem \ref{potmod}). 

\begin{ithm}\label{thma} Let $\{ r_\lambda\}$ be a weakly compatible system of $n$-dimensional $l$-adic representations of $G_F$ defined over $M$ and unramified outside $S$, where for simplicity we
assume  that $M$ contains the image of each embedding $F \into \barM$. Suppose that $\{ r_\lambda\}$ satisfies the following properties.
\begin{enumerate}
\item {\bf (Purity)} There is an integer $w$ such that, for each prime $v \not\in S$ of $F$, the roots of the common characteristic polynomial of the $r_\lambda(\Frob_v)$ are Weil $(\# k(v))^{w}$-numbers.
\item {\bf (Regularity)} For each embedding $\tau:F \into M$ the representation $r_\lambda$ has $n$ distinct $\tau$-Hodge--Tate numbers. 
\item {\bf (Odd essential self-duality)} $F$ is totally real; and either each $r_\lambda$ factors through a map to $\GSp_n(\barM_\lambda)$ with a totally odd multiplier character; or each $r_\lambda$ factors through a map to $\GO_n(\barM_\lambda)$ with a totally even multiplier character. Moreover in either case the multiplier characters 
form a weakly compatible system. 
\end{enumerate}

Then there is a finite, Galois, totally real extension $F'/F$ over which all the $r_\lambda$'s become
automorphic. In particular for any embedding $\imath:M \into \C$ the partial L-function $L^S(\imath \{ r_\lambda \},s)$ has meromorphic continuation to the whole complex plane and satisfies the expected functional equation.
\end{ithm}

A similar result is proved in the case that $F$ is CM. 

A similar theorem was proved in \cite{blggt} with an irreducibility assumption in place of our purity assumption. For compatible systems arising from geometry irreducibility can, in practice, be hard to check, but purity is often known thanks to Deligne's theorem. For instance theorem \ref{thma} has the following consequence. (See corollary \ref{coh}.)

\begin{icor} Suppose that $m \in \Z_{\geq 0}$, that $F$ is a totally real field and that $X/F$ is a smooth projective variety such that for all $\tau:F \into \C$ and all $i=0,...,m$ we have
\[ \dim H^{i,m-i}((X \times_{F,\tau} \C)(\C), \C) \leq 1. \]
Then $\{ H^m_{et}(X \times_F \barF, \Q_l) \}$ is a strictly compatible system of $l$-adic representations and the $L$-function $L(\{ H^m_{et}(X \times_F \barF, \Q_l) \}, s)$ has meromorphic continuation to the whole complex plane and satisfies the expected functional equation. 
\end{icor}

There is also a version of this corollary for motives, which is applicable more generally. (See corollary \ref{motpm}.)

\begin{icor} Suppose that $F$ is a totally real field, that $M$ is a number field, and that $X$ and $Y$ are pure motivated motives (in the sense of section 4 of \cite{and}) over $F$ with coefficients in $M$ such that the $l$-adic realizations of $X$ and $Y$ form weakly compatible systems of $l$-adic representations $\CH(X)$ and $\CH(Y)$. Suppose also that
\begin{enumerate}
\item ({\bf self-duality}) $X \cong X^\vee \otimes Y$
\item ({\bf regularity}) and for all $\tau: F \into \C$, all $\tau':M \into \C$, and all $i$ and $j$, we have
\[ \dim_\C H^{i,j}(X \times_{F,\tau} \C) \otimes_{M \otimes \C, \tau' \otimes 1} \C \leq 1. \]
\end{enumerate}

Then $\CH(X)$ is a strictly compatible system of $l$-adic representations and the $L$-function $L(\CH(X), s)$ has meromorphic continuation to the whole complex plane and satisfies the expected functional equation. 
\end{icor}

We use Katz's theory of rigid local systems to construct many examples of motivated motives to which this corollary applies. (See corollary \ref{rmpm} and  the discussion that follows it.)

Our approach to theorem \ref{thma} is to apply the potential automorphy result theorem 4.5.1 of \cite{blggt} to the irreducible constituents of the $r_\lambda$. The arguments of section 5 of \cite{blggt} can be applied as long as one can show that \vspace{2mm}

\noindent {\em there is a positive Dirichlet density set $\CL$ of rational primes such that for all $\lambda|l \in \CL$ the irreducible constituents of $r_\lambda$ are all odd, essentially self-dual.}\vspace{2mm}

\noindent In \cite{blggt} this was proved subject to the very strong condition that $\{ r_\lambda\}$ is extremely regular. However this does not apply in many settings. The main innovation of the present paper is to prove this assertion under the assumption that $\{ r_\lambda\}$ is pure and regular. (See lemma \ref{dual}.) We do this using a technique from \cite{patrikis} (see \cite[\S $6$, \S $16.3$]{patrikis}), which makes use of the CM nature of the field of coefficients.
Curiously our set $\CL$ will only contain primes $l$, such that all primes $\lambda$ above $l$ in a certain CM field (a souped-up field of definition for $\{ r_\lambda\}$) satisfy ${}^c\lambda=\lambda$.  

The same technique, combined with the arguments of section 5 of \cite{blggt} also allow us to prove the following irreducibility result for $l$-adic representations arising from a cuspidal automorphic representation. (See theorem \ref{irred}.)

\begin{ithm} Suppose that $F$ is a CM (or totally real)
field and that $\pi$ is a polarizable, regular algebraic, cuspidal automorphic representation of $GL_n(\A_F)$. Then there is a positive Dirichlet density set $\CL$ such that, for all $l\in \CL$ and all $\imath:\barQQ_l \iso \C$, the representation $r_{l,\imath}(\pi)$ is irreducible.\end{ithm} 

In this paper will use, often without comment, the notation and definitions of \cite{blggt}, particularly of sections 2.1 and 5.1 of that paper. For instance $\epsilon_l$ will denote the $l$-adic cyclotomic character and $V(m)$ for the Tate twist $V(\epsilon_l^m)$.  Also by a CM field we will mean a number field $F$ admitting an automorphism $c$, which coincides with complex conjugation for every embedding $F \into \C$. If $F$ is a CM field then $F^+$ will denote $F^{ \{ 1,c\} }$ the maximal totally real subfield of $F$. We have $[F:F^+]=1$ or $2$.
Finally if $\Gamma \supset \Delta$ are groups, if $\gamma \in \Gamma$ normalizes $\Delta$ and if $r$ is a representation of $\Delta$ we define a representation $r^\gamma$ of $\Delta$ by
\[ r^\gamma(\delta)=r(\gamma \delta \gamma^{-1}). \]
 \newpage

\section{Generalities on compatible systems of $l$-adic representations.}\label{sec:main result}

We will recall some elementary facts and definitions concerning compatible systems of $l$-adic representations.

First of all recall that if $F/ \Q_l$ is a finite extension, if $V$ is a finite dimensional $\barQQ_l$-vector space, if 
\[ r: G_F \lra GL(V) \]
is a continuous semi-simple representation, and if $\tau:F \into \barQQ_l$ is a continuous embedding, then we define a multiset
\[ \HT_\tau(r) \]
of integers by letting $i$ have multiplicity
\[ \dim_{\barQQ_l} (V(\epsilon_l^i) \otimes_{\tau,F} \widehat{\barF})^{G_F}, \]
where $\widehat{\barF}$ denotes the completion of the algebraic closure of $F$. Then $\# \HT_\tau(r) \leq \dim_{\barQQ_l}V$ with equality if $r$ is de Rham (or even, by definition, Hodge-Tate).
For $\sigma \in \Gal(\barQQ_l/\Q_l)$ we have
\[ \HT_{\sigma \tau}({}^\sigma r)=\HT_\tau (r). \]
In particular if the trace, $\tr r$, is valued in a finite extension $M/\Q_l$ (in $\barQQ_l$) and if $\sigma \in \Gal(\barQQ_l/M)$ then 
\[ \HT_{\sigma \tau}(r)=\HT_\tau(r). \]
Also if $\sigma$ is a continuous automorphism of $F$, then
\[ \HT_{\tau \circ \sigma}(r^\sigma) = \HT_\tau(r). \] 

Suppose now that $F$ is a number field, that $V$ is a finite dimensional $\barQQ_l$-vector space and that
\[ r: G_F \lra GL(V) \]
is a continuous representation which is de Rham at all primes above $l$. If $\tau:F \into \barQQ_l$ then we define
\[ \HT_\tau(r) = \HT_\tau(r|_{G_{F_{v(\tau)}}}) \]
where $v(\tau)$ is the prime of $F$ induced by $\tau$. 
Note that if $\sigma \in G_{\Q_l}$ then 
\[ \HT_{\sigma \tau}({}^\sigma r) = \HT_\tau (r), \]
and so, if $\tr r$ is valued in a closed subfield $M \subset \barQQ_l$ and if $\sigma \in G_M$, then
\[ \HT_{\sigma \tau} (r)=\HT_\tau (r). \]
If $\sigma$ is an automorphism of $F$ and if $v$ is the prime of $F$ determined by $\tau:F \into \barQQ_l$, then 
\[ \begin{array}{rcl} (V^{\sigma}(\epsilon_l^i) \otimes_{\tau,F_{v(\tau)}} \widehat{\barF}_{v(\tau)})^{G_{F_{v(\tau)}}} &=& (V(\epsilon_l^i) \otimes_{\tau,F_{v(\tau)}} \widehat{\barF}_{v(\tau)})^{G_{F_{\sigma v(\tau)}}} \\ &\liso& (V(\epsilon_l^i) \otimes_{\tau\sigma^{-1},F_{v(\tau \sigma^{-1})}} \widehat{\barF}_{v(\tau \sigma^{-1})})^{G_{F_{v(\tau \sigma^{-1})}}}, \end{array} \]
where, in the middle space, $G_{F_{\sigma v(\tau)}}$ acts on $\widehat{\barF}_v$ via the conjugation by $\sigma^{-1}$ map $G_{F_{\sigma v(\tau)}} \ra G_{F_{v(\tau)}}$, and where the second map is $1 \otimes \sigma$. Thus we have
\[ \HT_\tau(r^\sigma)=\HT_\tau(r^\sigma|_{G_{F_{v(\tau)}}})= \HT_{\tau \sigma^{-1}}(r|_{G_{F_{v(\tau \sigma^{-1})}}}) = \HT_{\tau \sigma^{-1}}(r). \] 

Now let $F$ denote a number field. 
As in \cite{blggt}, by a {\em rank $n$ weakly compatible system of $l$-adic representations} $\CR$
{\em of} $G_F$ {\em defined over} $M$ we shall mean a $5$-tuple 
\[ (M,S,\{ Q_v(X) \}, \{r_\lambda \}, \{H_\tau\} ) \]
where
\begin{enumerate}
\item $M$ is a number field;
\item $S$ is a finite set of primes of $F$;
\item for each  prime $v\not\in S$ of $F$, $Q_v(X)$ is a monic degree $n$
polynomial in $M[X]$;
\item for each prime $\lambda$ of $M$ (with residue characteristic $l$ say) 
\[r_\lambda:G_F \lra \GL_n(\barM_\lambda) \]
is a continuous, semi-simple, representation such that 
\begin{itemize}
\item if $v \not\in S$ and $v \ndiv l$ is a prime of $F$ then $r_\lambda$
is unramified at $v$ and $r_\lambda(\Frob_v)$ has characteristic
polynomial $Q_v(X)$,
\item while if $v|l$ then $r_\lambda|_{G_{F_v}}$ is de Rham  and in the case $v \not\in S$ crystalline;
\end{itemize}
\item for $\tau:F \into \barM$, $H_\tau$ is a multiset of $n$ integers such that 
for any $i:\barM \into \barM_\lambda$ over $M$ we have 
$\HT_{i \circ\tau}(r_\lambda)=H_\tau$.
\end{enumerate}
We refer to a rank $1$ weakly compatible system of $l$-adic representations as a weakly compatible system of $l$-adic characters. 

Note that if $(M,S,\{ Q_v(X) \}, \{r_\lambda \}, \{H_\tau\} )$ is a weakly compatible system of $l$-adic representations of $G_F$ and that if $M' \supset M$ is a finite extension then the tuple $(M',S,\{ Q_v(X) \}, \{r_\lambda \}, \{H_\tau\} )$ is also a weakly compatible system of $l$-adic representations of $G_F$. Also note if $(M,S,\{ Q_v(X) \}, \{r_\lambda \}, \{H_\tau\} )$ is a weakly compatible system of $l$-adic representations of $G_F$ and if $\sigma \in G_M$, then $H_{\sigma \tau}=H_{\tau}$. (As they both equal $\HT_{i \circ \sigma \circ \tau}(r_\lambda)$, where $i:\barM \into \barM_\lambda$ over $M$.) 

We will say that two weakly compatible systems of $l$-adic representations of $G_F$ over $M$, say 
\[ \CR=(M,S,\{ Q_v(X) \}, \{r_\lambda \}, \{H_\tau\} )\]
 and 
 \[ \CR'=(M,S',\{ Q_v'(X) \}, \{r'_\lambda \}, \{H_\tau'\} ),\]
  are {\em equivalent}  if $Q_v(X)=Q_v'(X)$ for a set of $v$ of Dirichlet density $1$. We write $\CR \equiv \CR'$. In this case we have $Q_v(X)=Q_v'(X)$ for all $v \not\in S \cup S'$, we have $r_\lambda'\cong r_\lambda$ for all $\lambda$, and we have $H_\tau=H_\tau'$ for all $\tau$. 

If $\CR=(M, S,\{ Q_{v}(X) \}, \{r_\lambda \}, \{H_\tau\} )$ is a weakly compatible system of $l$-adic representations of $G_F$ and if $\sigma \in \Aut(M)$ then we set 
\[ {}^\sigma \CR=(M,S,\{ {}^\sigma Q_v(X) \}, \{{}^\sigma r_{\sigma^{-1} \lambda} \}, \{H_{\sigma^{-1} \tau}\} ). \]
It is again a weakly compatible system of $l$-adic representations of $G_F$. Similarly if $\sigma \in \Aut(F)$ then we set
\[ \CR^\sigma=(M,\sigma^{-1}S,\{ Q_{\sigma v}(X) \}, \{r_{\lambda}^\sigma \}, \{H_{\tau \sigma^{-1}}\} ). \]
It is again a weakly compatible system of $l$-adic representations of $G_F$. 

\begin{lem} Suppose that $\CR=(M,S,\{ Q_v(X) \}, \{r_\lambda \}, \{H_\tau\} )$ is a weakly compatible system of $l$-adic representations of $G_F$ and that $M' \subset M$ is a subfield such that $M'[X]$ contains $Q_v(X)$ for a set of $v\not\in S$ of Dirichlet density $1$. Then $(M',S,\{ Q_v(X) \}, \{r_\lambda \}, \{H_\tau\} )$ is a weakly compatible system of $l$-adic representations of $G_F$. \end{lem}

\pfbegin Enlarging $M$ if need be we may assume that $M/\Q$ is Galois. If $\lambda'$ is a prime of $M'$ we define $r_{\lambda'}$ to be $r_\lambda$ for any prime $\lambda|\lambda'$ of $M$. To see that this is well defined we must check that if $\lambda_1,\lambda_2$ are primes of $M$ above $\lambda'$ and if $\sigma:\barM_{\lambda_1} \iso \barM_{\lambda_2}$ is an isomorphism restricting to the identity on $M'_{\lambda'}$ then ${}^\sigma r_{\lambda_1} \cong r_{\lambda_2}$. This follows because these two representations have the same trace, which fact in turn follows from the Cebotarev density theorem because $\tr {}^\sigma r_{\lambda_1}(\Frob_v) = \tr r_{\lambda_2}(\Frob_v)$ for a set of primes $v$ of $F$ of Dirichlet density $1$. 

We must next check that if $v \not\in S$ then $Q_v(X) \in M'[X]$. If $\lambda'$ is a prime of $M'$ and $\sigma \in G_{M'_{\lambda'}}$ then ${}^\sigma r_{\lambda'} \cong r_{\lambda'}$ so that ${}^\sigma Q_v(X)=Q_v(X)$. As $\Gal(M/M')$ is generated by the elements of decomposition groups at finite primes we see that $Q_v(X)$ is fixed by $\Gal(M/M')$ and the claim follows.

Finally we must check that if $\tau:F \into \barM$ and $i:\barM \into \barM_{\lambda'}'$  over $M'$ then $\HT_{i \circ \tau}(r_{\lambda'})=H_\tau$. Choose $\sigma \in G_{M'}$ such that $i \circ \sigma^{-1}$ is $M$-linear.  Then we certainly have
\[ \HT_{i \circ \tau}(r_{\lambda'}) = H_{\sigma \tau}. \]
Thus it suffices to check that for all $\sigma \in G_{M'}$ and all $\tau:F \into \barM$ we have
\[ H_{\sigma \tau} = H_\tau. \]
In fact we only need treat the case that $\sigma \in G_{M'_{\lambda_0'}}$ for some prime $\lambda'_0$. (As such elements topologically generate $G_{M'}$.)

Let $\lambda_0$ be the prime of $M$ above $\lambda'_0$ corresponding to a given embedding $G_{M'_{\lambda_0'}} \into G_{M'}$, and let $j:\barM \into \barM_{\lambda_0}$ be $M$-linear. Then
\[ r_{\lambda_0} \cong {}^\sigma r_{\sigma^{-1}\lambda_0} \]
and so
\[ H_{\sigma \circ \tau} = \HT_{j \circ \sigma \circ \tau} (r_{\lambda_0}) = \HT_{j \circ \sigma \circ \tau}({}^\sigma r_{\sigma^{-1}\lambda_0}) = \HT_{j' \circ \tau}(r_{\sigma^{-1} \lambda_0}) = H_\tau, \]
where
\[ j'=\sigma^{-1} \circ j \circ \sigma : \barM \into \barM_{\sigma^{-1} \lambda_0} \]
is $M$-linear.
\pfend

If the conclusion of this lemma holds we will say that $\CR$ {\em can be defined over $M'$}. The lemma implies that any weakly compatible system of $l$-adic representations has a unique minimal field of definition, namely the subfield of $M$ generated by the coefficients of all the $Q_v(X)$ for $v \not\in S$. If $M/\Q$ is Galois then it is also the fixed field of
\[ \{ \sigma \in \Gal(M/\Q): \,\, {}^\sigma \CR \equiv \CR\}. \]

We will call $\CR$ {\em pure} of weight $w$ if for each $v \not\in S$, for each root $\alpha$ of $Q_v(X)$ in $\barM$ and for each $\imath:\barM \into \C$
we have 
\[ | \imath \alpha |^2 = (\# k(v))^w. \]
Note that this definition is apparently slightly weaker than the definition given in \cite{blggt}, but the next lemma, which is essentially due to \cite{patrikis} (see \cite[$16.3.1$, $16.3.3$]{patrikis}), shows that the two definitions are actually equivalent.

\begin{lem}\label{pure} Suppose that $\CR=(M,S,\{ Q_v(X) \}, \{r_\lambda \}, \{H_\tau\} )$ is a pure weakly compatible system of $l$-adic representations of $G_F$ of weight $w$. 
\begin{enumerate}
\item If $c$ is the restriction to $M$ of any complex conjugation then ${}^c \CR \cong \CR^\vee \otimes  \{ \epsilon_l^{-w} \}$.
\item $\CR$ can be defined over a CM field.
\item If $c \in \Aut(\barM)$ denotes a complex conjugation then 
\[ H_{c \tau} = \{ w-h: \,\, h \in H_\tau \}. \]
\item Suppose that $F_0$ is the maximal CM subfield of $F$. If $\tau|_{F_0}=\tau'|_{F_0}$ then
\[ H_\tau = H_{\tau'}. \]
\end{enumerate}
\end{lem}

\pfbegin The first part follows as for $v \not\in S$ we have
\[ {}^c Q_v(X)=X^{\dim \CR} Q_v(q_v^w/X)/Q_v(0), \]
where $q_v=\# k(v)$. If $M/\Q$ is Galois and $c,c' \in \Gal(M/\Q)$ are two complex conjugations, we deduce that 
\[ {}^{cc'} \CR \cong \CR. \]
The second part follows.
The third part also follows from the first.

For the fourth part we may assume that $M$ is a CM field (using the second part), that $M/\Q$ is Galois and that $M$ contains $\tau F_0$ for all $\tau:F \into \barM$. Then for any $\tau:F \into \barM$ we see that $M \cap \tau F =\tau F_0$ and that $M(\tau F) \cong M \otimes_{\tau, F_0} F$. Thus if $\tau|_{F_0}=\tau'|_{F_0}$ we see that $M (\tau F) \cong M (\tau' F)$ as $M \otimes F$-algebras so that there is $\sigma \in G_M$ with $\sigma \tau = \tau'$. The fourth part follows.
\pfend 

Recall (from section 5.1 of \cite{blggt}) that one can apply standard linear algebra operations like direct sum, tensor product and dual to compatible systems of $l$-adic representations. One can also restrict them from $G_F$ to $G_{F'}$ if $F' \supset F$.

\begin{lem}\label{charht} Suppose that $F$ is a CM field, that $\CR=(M,S,\{ Q_v(X) \}, \{r_\lambda \}, \{H_\tau\} )$ is a pure weakly compatible system of $l$-adic representations of $G_F$ of weight $w$ and that $\CM=(M,S^\mu,\{X-\alpha_v\}, \{ \mu_\lambda \}, \{ H^\mu_\tau \})$ is a weakly compatible system of $l$-adic characters of $G_{F^+}$ with
\[ \CR^c \equiv \CM|_{G_F} \otimes \CR^\vee . \]
Then for all $\tau$ we have $H^\mu_\tau=\{ w\}$. \end{lem}

\pfbegin
This follows on noting that $\CM$ must be pure of weight $2w$ and using the classification of algebraic $l$-adic characters of the absolute Galois group of a totally real field. (They must all be of the form $\mu_0 \epsilon_l^{-w_0}$ where $\mu_0$ has finite order and $w_0 \in \Z$.) 
\pfend

Recall that we call $\CR$ {\em regular} if for each $\tau:F \into \barM$ every element of $H_\tau$ has multiplicity $1$. 

\begin{lem} Suppose that $\CR=(M,S,\{ Q_v(X) \}, \{r_\lambda \}, \{H_\tau\} )$ is a regular pure weakly compatible system of $l$-adic representations of $G_F$ of weight $w$. Then we may replace $M$ by a CM field $M'$ such that for all open subgroups $H$ of $G_F$ and all primes $\lambda$ of $M'$, all sub-representations of $r_\lambda|_H$ are defined over $M'_\lambda$. \end{lem}

\pfbegin By lemma \ref{pure} we may assume $M$ is CM. Then this lemma is proved in the same way as lemma 5.3.1(3) of \cite{blggt}, 
noting that, by purity, the splitting field over $M$ for the polynomial $Q_v(X)Q_{v'}(X)$ which occurs in the proof of lemma 5.3.1(3) of \cite{blggt} is a CM field.
\pfend

When the conclusion of this lemma holds for $M$ we will call $M$ a {\em full field of definition} for $\CR$.

If $M$ is a CM  field we will call a prime $\lambda$ of $M$ {\em conjugation invariant} if ${}^c\lambda=\lambda$. Thus if $M$ is totally real all primes of $M$ are conjugation invariant.

\begin{lem} Suppose that $M$ is a CM field. Then there is a set of rational primes $\Sigma$ of positive Dirichlet density such that all primes of $M$ above any element $l \in \Sigma$ are conjugation invariant. \end{lem}

\pfbegin One can reduce to the case that $M/\Q$ is Galois. In this case one can for instance take $\Sigma$ to be the set of rational primes $l$ unramified in $M$ and with $[\Frob_l]=\{c\} \subset \Gal(M/\Q)$. \pfend

\begin{lem}\label{dual} Suppose that $F$ is a CM field, that $\CR=(M,S,\{ Q_v(X) \}, \{r_\lambda \}, \{H_\tau\} )$ is a regular, pure weakly compatible system of $l$-adic representations of $G_F$ of weight $w$ and that $\CM=(M,S^\mu,\{X-\alpha_v\}, \{ \mu_\lambda \}, \{ H^\mu_\tau \})$ is a weakly compatible system of characters of $G_{F^+}$ with
\[ \CR^c \equiv \CM|_{G_F} \otimes \CR^\vee . \]
Suppose further that $M$ is a CM field and is a full field of definition for $\CR$.
Let $c$ denote a complex conjugation in $G_{F^+}$. Suppose further that $\lambda$ is a conjugation invariant prime of $M$ and that $\langle \,\,\,,\,\,\, \rangle$ is a bilinear form on the space underlying $r_\lambda$ such that 
\[ \langle r_\lambda(\sigma) x, r_\lambda(c \sigma c^{-1}) y\rangle = \mu_\lambda(\sigma) \langle x,y\rangle \]
for all $\sigma \in G_F$ and all $x,y$ in the underlying space of $r_\lambda$. 

Then the irreducible constituents of $r_\lambda$ are orthogonal with respect to $\langle \,\,\,,\,\,\,\rangle$.
\end{lem}

\pfbegin
Note that $c$ on $M$ extends to a unique continuous automorphism of $M_\lambda$. Let $r$ be a constituent of $r_{\lambda}$.
For $v \not\in S$ we let $Q^r_v(X)\in M_\lambda[X]$ denote the characteristic polynomial of $r(\Frob_v)$. Then by purity we see that
\[ {}^c Q^r_v(X)=X^{\dim r} Q_v^r(q_v^w/X)/Q_v(0) \]
where $q_v = \# k(v)$. Using the Cebotarev density theorem we deduce that
\[ {}^cr \cong \epsilon_l^{-w} r^\vee. \]
Hence
\[ \mu_\lambda (r^c)^\vee \cong (\mu_\lambda \epsilon_l^w)\,\, {}^cr^c. \]
Using lemma \ref{charht} we deduce that for $\tau:F \into \barM_\lambda$ we have
\[ \HT_\tau(\mu_\lambda (r^c)^\vee)= \HT_\tau((\mu_\lambda \epsilon_l^w)\,\, {}^cr^c) = \HT_\tau({}^cr^c) = \HT_{c\tau c}(r) = \HT_\tau(r). \]
By regularity we deduce that
\[ r \cong   \mu_\lambda (r^c)^\vee, \]
and moreover $r_\lambda$ has no other sub-representation isomorphic to $\mu_\lambda (r^c)^\vee$. The lemma follows.
\pfend

We remark that the lemma is presumably true without the assumption that $\lambda$ is conjugation invariant, but we don't know how to prove this.

We will call $\CR$ {\em strictly compatible} if for each finite place $v$ of $F$ there is a Weil--Deligne representation $\WD_v(\CR)$ of $W_{F_v}$ over $\barM$ such that for each place $\lambda$ of $M$ 
and every $M$-linear embedding $\varsigma:\barM \into \barM_\lambda$ the push forward $\varsigma\WD_v(\CR)\cong \WD(r_\lambda|_{G_{F_v}})^\Fsemis$. (This is slightly stronger than the notion we defined in \cite{blggt}.) Moreover we will call $\CR$ {\em strictly pure} of weight $w$ if
 $\CR$ is strictly compatible and for each prime $v$ of $F$ the Weil--Deligne representation $\WD_v(\CR)$ is pure of weight $w$ (see section 1.3 of \cite{blggt}).
 
 If $\CR$ is pure and if $\imath:M \into \C$ then we can define the partial $L$-function $L^S(\imath \CR,s)$, which is a holomorphic function in some right half-plane. If $\CR$ is regular and strictly pure we can the completed $L$-function $\Lambda(\imath \CR,s)$ and the epsilon factor $\epsilon(\imath \CR,s)$. (See section 5.1 of \cite{blggt} for details.)
 
 Suppose that $F$ is a number field and that $\pi$ is a regular algebraic cuspidal automorphic representation of $GL_n(\A_F)$ of weight $a \in (\Z^n)^{\Hom(F,\C),+}$. 
(See section 2.1 of \cite{blggt}.) We define $M_\pi$ to be the fixed field of the subgroup of $\Aut(\C)$ consisting of elements $\sigma \in \Aut(\C)$ such that
${}^\sigma \pi^\infty \cong \pi^\infty$. It follows from theorem 3.13 of \cite{clozel} that $M_\pi$ is a number field, and in fact is a CM field \cite[$6.2.3$]{patrikis}. We let $S_\pi$ denote the set of primes of $F$ where $\pi$  is ramified and for $v \not\in S_\pi$ a prime of $F$ we let $Q_{\pi,v}(X)\in M_\pi[X]$ denote the characteristic polynomial of $\rec_{F_v}(\pi_v |\det|_v^{(1-n)/2})$. If $\barM_\pi$ denotes the algebraic closure of $M$ in $\C$ and if $\tau: F \into \barM_\pi$ then we set
\[ H_{\pi,\tau}=\{ a_{\tau,1}+n-1,a_{\tau,2}+n-2,...,a_{\tau,n}\}. \]
In the case that $F$ is a CM field and $\pi$ is polarizable, it is known that for each prime $\lambda$ of $M_\pi$ there is a continuous semi-simple representation
\[ r_{\pi,\lambda}: G_F \lra GL_n(\barM_{\pi,\lambda}) \]
such that
\[ \CR_\pi=(M_\pi, S_\pi, \{ Q_{\pi,v}(X) \}, \{ r_{\pi,\lambda} \}, \{ H_\tau \}) \]
is a regular, strictly pure compatible system of $l$-adic representations. (Combine theorem 2.1.1 of \cite{blggt}, theorem 1.1 of \cite{ana} and the usual twisting and descent arguments. Note that theorem 2.1.1 of \cite{blggt} simply collects together results of many other authors, see that paper for more details.)

Recall from section 5.1 of \cite{blggt} that a weakly compatible system of $l$-adic representations $\CR=(M,S,\{ Q_v(X) \}, \{r_\lambda \}, \{H_\tau\} )$ is called {\em automorphic} if there is a regular algebraic, cuspidal automorphic representation $\pi$
of $\GL_n(\A_F)$ and an embedding $\imath:M \into \C$, such that if $v \not\in S$ then $\pi_v$ is unramified and 
$\rec(\pi_v |\det|_v^{(1-n)/2})(\Frob_v)$ has characteristic
polynomial $\imath (Q_v(X))$. 

Suppose that $F$ is a CM field, that $\CR=(M,S,\{ Q_v(X) \}, \{r_\lambda \}, \{H_\tau\} )$ is a weakly compatible system of $l$-adic representations of $G_F$ and that $\CM=(M,S^\mu,\{X-\alpha_v\}, \{ \mu_\lambda \}, \{ H^\mu_\tau \})$ is a weakly compatible system of $l$-adic characters of $G_{F^+}$. Recall (from section 5.1 of\cite{blggt}) that we 
call $(\CR,\CM)$ a {\em polarized} (resp. {\em totally odd, polarized}) weakly compatible system if for all primes $\lambda$ of $M$ the pair $(r_\lambda,\mu_\lambda)$ is a polarized (resp. totally odd polarized) $l$-adic representation in the sense of section 2.1 of \cite{blggt}. 

If $F$ is a CM field and $(\CR,\CM)$ is an automorphic, regular, polarized weakly compatible system of $l$-adic representations of dimension $n$, then $\CR\equiv \CR_\pi$ for some regular algebraic, polarizable, cuspidal automorphic representation $\pi$ of $GL_n(\A_F)$.

We finish this section with an application of lemma \ref{dual} to irreducibility results.

\begin{thm} \label{irred} Suppose that $F$ is a CM  
field and that $\pi$ is a polarizable, regular algebraic, cuspidal automorphic representation of $GL_n(\A_F)$. Then there is a finite CM extension $M/M_\pi$ and a Dirichlet density $1$ set $\CL$ of rational primes, such that for all conjugation-invariant primes $\lambda$ of $M$ dividing an $\ell \in \CL$, $r_{\pi,\lambda|_{M_\pi}}$ is irreducible. 

In particular, there is a positive Dirichlet density set $\CL'$ of rational primes such that if a prime $\lambda$ of $M_{\pi}$ divides some $\ell \in \CL'$, then $r_{\pi, \lambda}$ is irreducible.
\end{thm} 

\pfbegin
The proof is the same as the proof of theorem 5.5.2 and proposition 5.4.6 of \cite{blggt}, except that instead of appealing to lemma 5.4.5 of \cite{blggt} one appeals to lemma \ref{dual} above.
\pfend

\newpage

\section{Applications}

We will first consider applications of lemma \ref{dual} to potential automorphy theorems.

\begin{thm}\label{potmod} Suppose that $F/F_0$ is a finite Galois extension of CM 
  fields, and that $\Favoid/F$ is a finite Galois extension. Suppose also for $i=1,...,r$ that $(\CR_i,\CM_i)$ is a totally odd, polarized weakly compatible system of $l$-adic representations with each $\CR_i$ pure and regular. Then there is a finite CM  extension $F'/F$ such that $F'/F_0$ is Galois and $F'$ is linearly disjoint from $\Favoid$ over $F$, with the following property. For each $i$ we have a decomposition 
  \[ \CR_i = \CR_{i,1} \oplus ... \oplus \CR_{i,s_i} \]
  into weakly compatible systems $\CR_{i, j}$ with each 
  $(\CR_{i,j},\CM_i)|_{G_{F'}}$ automorphic.
  \end{thm}
   
\pfbegin
Suppose that
\[ \CR_i=(M_i,S_i,\{ Q_{i,v}(X) \}, \{r_{i,\lambda} \}, \{H_{i,\tau}\} ) \]
and that $M_i$ is a CM field which is a full field of definition for $\CR_i$.
Write
\[ r_{i,\lambda} = \bigoplus_j r_{i,\lambda,j} \]
with each $r_{i,\lambda,j}$ irreducible. It follows from lemma \ref{dual} that if $\lambda$ is a conjugation invariant prime of $M_i$ then $(r_{i,\lambda,j},\mu_{i,\lambda})$ is a totally odd, polarized $l$-adic representation.

Let $\CL$ be the Dirichlet density $1$ set of rational primes obtained by applying proposition 5.3.2 of \cite{blggt} to $\CR$. Then for $\lambda|l \in \CL$ we see that
$\barr_{i,\lambda,j}|_{G_{F(\zeta_l)}}$ is irreducible. 
Removing finitely many primes from $\CL$ we may further suppose that 
\begin{itemize}
\item $l \in \CL$ implies $l\geq 2(\dim \CR_i+1)$ for all $i$,
\item $l\in \CL$ implies $l$ is unramified in $F$ and $l$ lies below no element of any $S_i$,
\item if $\lambda|l \in \CL$ then all the Hodge--Tate numbers of each $r_{i,\lambda}$ lie in a range of the form $[a,a+l-2]$.
\end{itemize}
From proposition 2.1.2 of \cite{blggt}, we deduce that for all $\lambda|l \in \CL$, the image $\barr_{i,\lambda,j}(G_{F(\zeta_l)})$ is adequate. Moreover, by lemma 1.4.3 of \cite{blggt}, $r_{i,\lambda,j}$ is potentially diagonalizable. 

For each $i=1,...,r$ choose a prime $l_i \in \cL$ and a conjugation invariant prime $\lambda_i|l_i$ of $M_i$. Replace $\Favoid$ by its compositum with the $\barF^{\ker \barr_{i,\lambda_i,j}}(\zeta_{l_i})$ for all $i,j$. Now apply theorem 4.5.1 of \cite{blggt} to $\{ r_{i,\lambda_i,j} \}$. We obtain a CM extension $F'/F$ with $F'/F_0$ Galois and $F'$ linearly disjoint from $\Favoid$ over $F$; and regular algebraic, cuspidal, polarized automorphic representations $(\pi_{i,j},\chi_{i,j})$ of $GL_{n_{i,j}}(\A_{F'})$ such that
\[ r_{l_i,\imath_i}(\pi_{i,j}) = r_{i,\lambda_i,j}|_{G_{F'}} \]
and
\[r_{l,\imath_i}(\chi_i) \epsilon_{l_i}^{1-n_{i,j}} = \mu_{i,\lambda_i} \]
for some $\imath_i:\barM_{i,\lambda_i} \iso \C$. 

As in the proof of Theorem 5.5.1 of \cite{blggt} we see that $r_{i,\lambda_i,j}$ is part of a weakly compatible system $\CR_{i,j}$ of $l$-adic representations of $G_F$. It is immediate that 
\[ \CR_i \cong \bigoplus_j \CR_{i,j} \]
and that each $(\CR_{i,j},\CM_i)|_{G_{F'}}$ is automorphic.
\pfend

The next corollary now follows in the same way as corollary 5.4.3 of \cite{blggt}.

\begin{cor}\label{compatiblesystem}  Suppose that $F$ is a CM 
  field, and that $(\CR,\CM)$ is a totally odd, polarized weakly compatible system of $l$-adic representations of $G_F$ with $\CR$ pure and regular. 
  \begin{enumerate}
\item If $\imath:M \into \C$, then $L^S(\imath \CR,s)$ converges
  (uniformly absolutely on compact subsets) on some right half plane
  and has meromorphic continuation to the whole complex plane.
\item The compatible system $\CR$ is strictly pure. Moreover
\[ \Lambda(\imath \CR,s)=\epsilon(\imath \CR,s) \Lambda(\imath \CR^\vee, 1-s). \]
\end{enumerate} \end{cor}

\begin{cor}\label{coh} Suppose that $m \in \Z_{\geq 0}$, that $F$ is a totally real field and that $X/F$ is a smooth projective variety such that for all $\tau:F \into \C$ and all $i=0,...,m$ we have
\[ \dim H^{i,m-i}((X \times_{F,\tau} \C)(\C), \C) \leq 1. \]
Then $\{ H^m_{et}(X_{/ \barF}, \Q_l) \}$ is a strictly pure compatible system of $l$-adic representations and the $L$-function $\Lambda(\{ H^m_{et}(X_{/\barF}, \Q_l) \}, s)$ has meromorphic continuation to the whole complex plane and satisfies the functional equation
\[ \Lambda(\{ H^m_{et}(X_{/\barF}, \Q_l) \},s)=\epsilon(\{ H^m_{et}(X_{/ \barF}, \Q_l) \},s) \Lambda(\{ H^m_{et}(X _{/\barF}, \Q_l) \}, 1+m-s). \]
\end{cor}

\pfbegin
By hard Lefschetz and Poincar\'{e} duality 
we see that there is a perfect pairing 
\[ H^m_{et}(X_{/ \barF}, \Q_l) \times H^m_{et}(X_{/\barF}, \Q_l) \lra \Q_l(-m) \]
of parity $(-1)^m$. 
Thus 
\[ ( \{ H^m_{et}(X_{/\barF}, \Q_l) \}, \{ \epsilon_l^{-m} \}) \]
is a totally odd, polarized weakly compatible system of $l$-adic representations of $G_F$. Moreover $\{ H^m_{et}(X_{/\barF}, \Q_l) \}$ is pure, by Deligne (\cite{deligne}), and regular. This corollary follows from the previous one and the isomorphism
\[ H^m_{et}(X_{/\barF}, \Q_l)^\vee \cong H^m_{et}(X_{/\barF}, \Q_l)(m). \]
\pfend

We can extend this corollary to motives. We choose to do this in Andr\'{e}'s category of motivated motives. 
Suppose that $F$ and $M$ are number fields. We will let $\CM_{F,M}$ denote the category of motivated motives over $F$ with coefficients in $M$ (see section 4 of \cite{and}). An object of $\CM_{F,M}$ will be a triple $(X,p,m)$, where 
\begin{itemize}
\item $X/F$ is a smooth, projective variety,
\item $p \in C^0_{mot}(X,X)_M$ is an idempotent,
\item and $m$ is an integer.
\end{itemize}
(See section 4 of \cite{and}, and section 2 of \cite{and} for the definition of $C^0_{mot}(X,X)_M$.) 
If $(X,p,m)$ is an object of $\CM_{F,M}$,
then we can form the following cohomology groups:
\begin{itemize} 
\item The de Rham realization
\[
H_\DR(X,p,m) = \bigoplus_{i} H^i_\DR(p) \left( H^i_\DR(X/F) \otimes_{\Q} M \right),
\]
an $(F \otimes_\Q M)$-module with an exhaustive, separated filtration 
\[ \Fil^k(X,p,m) = \bigoplus_{i} H^i_\DR(p) \Fil^{k+m} \left( H^i_\DR(X/F) \otimes_{\Q} M \right) \]
 by sub-$(F \otimes_\Q M)$-modules.
\item For all $\tau \colon F \into \C$, the $\tau$-Betti realization 
\[
H_{B,\tau}(X,p,m)=\bigoplus_{i} H^i_{B,\tau}(p)H^i_{B}\left( (X \times_{F, \tau} \C)^{an}, M \right)(m),
\] 
an $M$-vector space. 
\item The $\lambda$-adic realization 
\[
H_\lambda(X,p,m) = \bigoplus_{i} H_\lambda^i(p)H^i_{et}(X \times_F \bar{F}, M_{\lambda})(m),
\] 
an $M_\lambda$-vector space with a continuous $G_F$-action.
\end{itemize}
There is also a natural $M$-linear map
\[ c: H_{B,\tau}(X,p,m) \lra H_{B,c\tau}(X,p,m). \]
If $\ttau:\barF \into \C$ extends $\tau$ then we get a comparison isomorphism
\[ \alpha_{\lambda,B, \ttau}: H_\lambda(X,p,m) \liso H_{B,\tau}(X,p,m) \otimes_M M_\lambda. \]
For $\sigma \in G_F$ we have $\alpha_{\lambda,B,\ttau \sigma}=\alpha_{\lambda,B,\ttau} \circ \sigma$. We also have $\alpha_{\lambda,B,c \ttau} = c \circ \alpha_{\lambda,B,\ttau}$. We also have a comparison map
\[ \alpha_{\DR,B,\tau}: H_\DR(X,p,m) \otimes_{F,\tau} \C \liso H_{B,\tau}(X,p,m) \otimes_\Q \C, \]
which satisfies
\[ (c \otimes c) \circ \alpha_{\DR,B,c\tau} = \alpha_{\DR,B,\tau} \circ (1 \otimes c): H_\DR(X,p,m)\otimes_{F,c\tau} \C \liso H_{B,\tau}(X,p,m) \otimes_\Q \C. \]
We will write 
\[ \Fil^i_\tau(X,p,m)=\alpha_{\DR,B,\tau} (\Fil^i(X,p,m) \otimes_{F,\tau} \C), \]
a $M \otimes_\Q \C$-submodule of $H_B(X,p,m) \otimes_\Q \C$, so that
\[ (c \otimes c) \Fil^i_\tau(X,p,m) = \Fil^i_{c \tau} (X,p,m). \]
We will refer to $\dim_M H_{B,\tau}(X,p,m)$, which is independent of $\tau$, as the rank of $(X,p,m)$. 

As described in section 4 of \cite{and} one can form the dual $X^\vee$ of an object $X$ of $\CM_{F,M}$, and the direct sum $X \oplus Y$ and the tensor product $X \otimes Y$ of two objects $X,Y$ of $\CM_{F,M}$. This makes $\cM_{F,M}$ a Tannakian category. The functors $H_\DR$, $H_{B,\tau}$ and $H_\lambda$ are faithful and exact. 
If $F'/F$ is a finite extension, one can restrict an object $X$ of $\CM_{F,M}$ to an object $X|_{F'}$ of $\CM_{F',M}$. 

We will call $(X,p,m)$ {\em compatible} if there is a finite set of primes $S$ of $F$ and, for each $v \not\in S$, a polynomial $Q_v\in M[X]$ such that: If $v \not\in S$ and $\lambda$ does not divide the residue characteristic of $v$, then $H_\lambda(X,p,m)$ is unramified at $v$ and $Q_v$ is the characteristic polynomial of $\Frob_v$ on $H_\lambda(X,p,m)$. If $(X,p,m)$ is compatible then
\[ \CH(X,p,m)=(M,S,\{ Q_v\}, \{H_\lambda(X,p,m)\}, \{ H_\tau \} ) \]
is a weakly compatible system of $l$-adic representations, where $H_\tau$ contains $i$ with multiplicity
\[ \dim_{\barM} \gr^i H_\DR(X,p,m) \otimes_{F \otimes M, \tau \otimes 1} \barM. \]
(To see that $\HT_\tau(H_\lambda(X,p,m))=H_\tau$ one can apply the remarks of section 2.4 of \cite{and}.)

We will call $(X,p,m)$ {\em pure of weight $w$} if (for instance) 
\[ H_\DR(X,p,m)= H_\DR^{w+2m}(p) \left( H_\DR^{w+2m}(X/F) \otimes_{\Q} M \right).\]
 Any object $(X,p,m)$ of $\CM_{F,M}$ can be written uniquely as
 \[ (X,p,m) = \bigoplus_r \mathrm{Gr}_r^W (X,p,m), \]
where $\mathrm{Gr}_r^W (X,p,m)$ is pure of weight $r$. (See section 4.4 of \cite{and}.)
The motivated motive $(X,p,m)$ is pure of weight $w$ if and only if $H_\lambda(X,p,m)$ is. In particular, if $(X,p,m)$ is compatible and pure, then $\CH(X,p,m)$ is also pure of the same weight.  If $(X,p,m)$ is pure of weight $w$ we have
\[ \begin{array}{rcl} H_{B,\tau}(X,p,m) \otimes_\Q \C &=& \Fil^i_\tau(X,p,m) \oplus (1 \otimes c) \Fil^{w+1-i}_\tau(X,p,m) \\ &=& \Fil^i_\tau(X,p,m) \oplus (c \otimes 1) \Fil^{w+1-i}_{c\tau}(X,p,m) . \end{array} \]
In particular $\gr^i_\tau(X,p,m)$ and $\gr^{w-i}_\tau(X,p,m) \otimes_{\C, c} \C$ are non-canonically isomorphic as $(M \otimes_\Q \C)$-modules.
Note that if $(X,p,m)$ is pure of weight $w$ and has rank $1$, and if $\tau:F \into \R$ then $w$ must be even and $\Fil_\tau^{w/2}(X,p,m) \neq \Fil^{w/2+1}_\tau(X,p,m)$.

We will call $(X,p,m)$ {\em regular} if for each embedding $\tau:F \into \barM$ and each $i \in \Z$ we have
\[ \dim_{\barM} \gr^i \left( H_\DR(X,p,m) \otimes_{F \otimes_\Q M, \tau \otimes 1_M} \barM \right) \leq 1. \]
If $(X,p,m)$ is regular and compatible, then $\CH(X,p,m)$ is regular.

The group $\Aut(F)$ also acts on the category $\CM_{F,M}$. An element $\sigma \in \Aut(F)$ takes an object $X$ of $\CM_{F,M}$ to ${}^\sigma X$. 
 We have the following observations:
 \begin{itemize}
 \item There is a $\sigma$-linear isomorphism $\sigma:H_\DR(X) \iso H_\DR({}^\sigma X)$ with $\sigma \Fil^i(X)=\Fil^i({}^\sigma X)$ for all $i$. 
 \item There is an isomorphism $\sigma: H_{B,\tau}(X) \iso H_{B,\tau \sigma^{-1}}({}^\sigma X)$.
 \item If we choose $\tsigma \in \Aut(\barF)$ lifting $\sigma$ then we get an isomorphism 
 \[ \tsigma: H_\lambda(X) \liso H_\lambda({}^\sigma X)\]
 such that $\tsigma \circ ( \tsigma^{-1}\sigma_1 \tsigma) =  \sigma_1  \circ \tsigma$ for all $\sigma_1 \in G_F$.
 \item We have
 \[ \alpha_{\lambda,B,\ttau \tsigma^{-1}} ({}^\sigma X)\circ \tsigma = \sigma \circ \alpha_{\lambda,B,\ttau}(X): H_\lambda(X) \lra H_{B,\tau \sigma^{-1}} ({}^\sigma X) \otimes_M M_\lambda , \]
 and
 \[ \alpha_{\DR,B,\tau \sigma^{-1}}({}^\sigma X) \circ \sigma = \sigma \circ \alpha_{\DR,B,\tau}(X): H_\DR(X) \otimes_{F,\tau} \C \lra H_{B, \tau \sigma^{-1}}({}^\sigma X) \otimes_\Q \C. \]
 \item If $X$ is compatible then $\CH({}^\sigma X)\cong \CH(X)^{\sigma^{-1}}$.
  \end{itemize}

If $F$ is a CM field, we will call an object $X$ of $\CM_{F,M}$ {\em polarizable} if there is an object $Y$ of $\CM_{F^+,M}$ (necessarily of rank $1$) such that
\[ X \cong Y|_F \otimes {}^cX^{\vee}. \]
We warn the reader that this is a non-standard use of the term `polarizable' in the context of the theory of motives. If $X$ is polarizable and $Y$ is as above, then 
there are non-degenerate pairings
\[ \langle \,\,\,,\,\,\,\rangle_\DR: H_\DR(X) \times H_\DR(X) \stackrel{1 \times c}{\lra} H_\DR(X) \times H_\DR({}^cX) \lra H_\DR(Y), \]
and
\[ \langle \,\,\,,\,\,\,\rangle_{B,\tau}: H_{B,\tau}(X) \times H_{B,\tau}(X) \stackrel{1 \times d}{\lra} H_{B,\tau}(X) \times H_{B,\tau}({}^cX) \lra H_{B,\tau}(Y), \]
where $d$ denotes the composite
\[ H_{B,\tau}(X) \stackrel{c}{\lra} H_{B,c\tau}(X)=H_{B,\tau c}(X) \stackrel{c}{\lra} H_{B,\tau}({}^c X). \]
We have
\[ \langle \alpha_{\DR,B,\tau}(X) x, (1 \otimes c)\alpha_{\DR,B,\tau}(X) y \rangle_{B,\tau}=\alpha_{\DR,B,\tau}(Y) \langle x,(1 \otimes c) y \rangle_\DR \]
for all $x,y \in H_\DR(X) \otimes_{F,\tau} \C$. 
Thus if $Y$ has weight $2w$ and $i+j>w$ then 
\[ \langle \Fil^i_\tau(X), (1 \otimes c) \Fil^j_\tau(X)\rangle_{B,\tau} =(0). \]
(Because, as $F^+$ is totally real and $Y$ has rank $1$, we have $\Fil^{w+1}(Y)=(0)$.)
Moreover if $\tilde{c} \in G_{F^+}$ lifts $c \in \Gal(F/F^+)$ then there is a non-degenerate pairing
\[ \langle \,\,\,,\,\,\,\rangle_{\lambda,\tilde{c}}: H_\lambda(X) \times H_\lambda(X) \stackrel{1 \times \tilde{c}}{\lra} H_\lambda(X) \times H_\lambda({}^cX) \lra H_\lambda(Y) \]
such that for all $\sigma \in G_F$ and $x,y \in H_\lambda(X)$ we have 
\[ \langle \sigma x , (\tilde{c}^{-1} \sigma \tilde{c}) y\rangle_{\lambda,\tilde{c}} = \sigma \langle x,y\rangle_{\lambda,\tilde{c}}. \]
If $\tau:F \into \C$ and $c(\tau)$ is the corresponding complex conjugation in $G_{F^+}$ we have
\[ \alpha_{\lambda,B,\tau} \langle x,y\rangle_{\lambda,c(\tau)} = \langle \alpha_{\lambda,B,\tau}x, \alpha_{\lambda,B,\tau} y \rangle_{B,\tau} \]
for all $x,y \in H_\lambda(X)$. 

\begin{lem} Suppose that $F$ is a CM field, that $M$ is a number field and that $X$ is a compatible, regular, polarizable, pure object of $\CM_{F,M}$. Then the weakly compatible system of $l$-adic representations $\CH(X)$ is totally odd, polarizable. \end{lem}

\pfbegin
It suffices to show that $\langle \,\,\,,\,\,\,\rangle_{\lambda,\tc}$ is symmetric for all complex conjugations $\tc \in G_{F^+}$. Choose an embedding $\tau_1:\barF \into \C$ such that $\tau_1 \circ \tc = c \circ \tau_1$. Then $\tau_1$ gives rise to an isomorphism
\[ H_{B,\tau_1|_F}(X) \otimes_M M_\lambda \cong H_\lambda(X). \]
Under this isomorphism $\langle \,\,\,,\,\,\,\rangle_{\lambda,\tc}$ correponds to $\langle \,\,\,,\,\,\,\rangle_{B,\tau_1|_F}$, so it suffices to show that $\langle \,\,\,,\,\,\,\rangle_{B,\tau}$ is symmetric for all $\tau:F \into \C$.

If $w$ denotes the weight of $X$ then we have
\[ \Fil^i_\tau(X) \oplus (1 \otimes c) \Fil^{w+1-i}_{\tau}(X) = H_{B,\tau}(X) \otimes_\Q \C. \]
Moreover as $Y$ must have weight $2w$, we see that, if $i+j>w$ then $\Fil^i_\tau(X)$ and $(1 \otimes c)\Fil^j_{\tau}(X)$ must annihilate each other under $\langle \,\,\,,\,\,\,\rangle_{B,\tau}$. 

Let $\tau':M \into \C$ and set
\[ H_{B,\tau}(X)_{\tau'}= H_{B,\tau}(X) \otimes_{M,\tau'} \C \]
and
\[ \Fil^i_\tau(X)_{\tau'}= \Fil^i_\tau(X) \otimes_{M \otimes_\Q\C, \tau' \otimes 1} \C. \]
Thus 
\[ \Fil^i_\tau(X)_{\tau'} \oplus (1 \otimes c) \Fil^{w+1-i}_{\tau}(X)_{c\tau'} = H_{B,\tau}(X)_{\tau'}. \]
Moreover $\langle \,\,\,,\,\,\,\rangle_{B,\tau}$ gives rise to a pairing
\[ \langle \,\,\,,\,\,\,\rangle_{B,\tau,\tau'}: H_{B,\tau}(X)_{\tau'}  \times H_{B,\tau}(X)_{\tau'} \lra H_{B,\tau}(Y)_{\tau'}, \]
under which $\Fil^i_\tau(X)_{\tau'}$ and $(1 \otimes c) \Fil^j_{\tau}(X)_{c \tau'}$ annihilate each other whenever $i+j>w$. It suffices to show that for all $\tau$ and $\tau'$ the pairing $\langle \,\,\,,\,\,\,\rangle_{B,\tau,\tau'}$ is symmetric.

Note that if $\gr^i_\tau(X)_{\tau'} \neq (0)$ then $\gr^{w-i}_{\tau}(X)_{c \tau'} \neq (0)$. Thus, by regularity,
\[ \Fil^i_\tau (X)_{\tau'} \cap (1 \otimes c) \Fil^{w-i}_{\tau}(X)_{c \tau'} \]
is one dimensional over $\C$. Let $e_i$ be a basis vector. As
\[ \Fil^{i+1}_\tau (X)_{\tau'} \cap (1 \otimes c) \Fil^{w-i}_{ \tau}(X)_{c \tau'}=(0) \]
we see that $e_i \in \Fil^i_\tau(X)_{\tau'}-\Fil^{i+1}_\tau(X)_{\tau'}$. Thus $\{ e_i \}$ is a basis of $H_{B,\tau}(X)_{\tau'}$. If $i>j$, then, since $e_i \in \Fil^i_\tau(X)_{\tau'}$ and $e_j \in (1 \otimes c) \Fil^{w-j}_{\tau}(X)_{c\tau'}$, 
\[ 
\langle e_i,e_j \rangle_{B,\tau,\tau'} = 0. 
\]
Similarly if $j>i$, then, since $e_i \in (1 \otimes c) \Fil^{w-i}_{ \tau}(X)_{c \tau'}$ and $e_j \in \Fil^j_\tau(X)_{\tau'}$, 
\[ \langle e_i,e_j \rangle_{B,\tau,\tau'} = 0. \]
We conclude that the matrix of $\langle \,\,\,,\,\,\,\rangle_{B,\tau,\tau'}$ with respect to the basis $\{ e_i\}$ is diagonal, and hence $\langle \,\,\,,\,\,\,\rangle_{B,\tau,\tau'}$ is symmetric, as desired.
\pfend

Combining this lemma with theorem \ref{potmod} we obtain the following corollary. 

\begin{cor}\label{motpm} Suppose that $F$ is a CM field, that $M$ is a number field and that $X$ is a compatible, regular, polarizable, pure object of $\CM_{F,M}$. 
Then there is a finite Galois CM extension $F'/F$ and a decomposition (perhaps after extending the field $M$)
\[ \CH(X) \equiv \CR_1 \oplus ... \oplus \CR_s \]
such that each $\CR_i|_{G_{F'}}$ is automorphic. 

In particular $\CH(X)$ is strictly pure and, if $\imath:M \into \C$, then $\Lambda(\imath \CH(X),s)$ has meromorphic continuation to the whole complex plane and satisfies the functional equation
\[ \Lambda(\imath \CH(X),s)=\epsilon(\imath \CH(X),s) \Lambda(\imath \CH(X)^\vee, 1-s). \]
\end{cor}

For example this corollary applies to the weight $n-1$ part of the motivated motive
$(Z_t,(1/\# H') \sum_{h \in H'} h,0)$, where
\begin{itemize}
\item $t \in \Q$;
\item $Z_t$ is the smooth hypersurface defined by \[X_0^{n+1}+ \ldots +X_n^{n+1} = (n+1)tX_0 \cdots X_n\] in projective $n$-space;
\item $H'$ is the 
group $\ker(\mu_{n+1}^{n+1} \stackrel{\prod}{\lra} \mu_{n+1})$, which acts on $Z_t$ by multiplication on the coordinates.
\end{itemize}
(In the case $t \not\in \Z[1/(n+1)]$ this was already proved in \cite{hsbt}.)

We next give a much more general example coming from Katz's theory of rigid local systems (\cite{katz:rls}). Let $S \subset \PP^1(\Q)$ be a finite set with complement $U=\PP^1-S$. Also let $N \in \Z_{>0}$. By a {\em rigid local system $\cF$ on $U_{\barQQ}$ with quasi-unipotent monodromy of order dividing $N$} we shall mean a lisse $\barQQ_l$-sheaf $\cF$ on $U_{\barQQ}$ with the following properties:
\begin{itemize}
\item $\cF$ is irreducible;
\item the $N^{th}$-power of the monodromy of $\cF$ at every point $s \in S$ is unipotent;
\item and any other $\barQQ_l$-local system on the complex manifold $U(\C)$ with monodromy conjugate to that of $\cF^\ann$ at every point $s \in S$ is isomorphimic to $\cF^\ann$.
\end{itemize}
By theorem 1.1.2 of \cite{katz:rls} the third condition is equivalent to the `cohomological rigidity' of $\cF$ in the sense of section 5.0 of \cite{katz:rls}. An abundant supply of such sheaves is supplied by the constructions of section 5.1 of \cite{katz:rls}. One may keep track of the monodromy at points $s \in S$ using the results of chapter 6 of \cite{katz:rls}. (See section 1.2 of \cite{dettweiler-reiter:G2rls} for a nice summary of these formulae.) 

Suppose that $\cF$ is a rigid local system with quasi-unipotent monodromy of order dividing $N$. According to theorem 8.4.1 of \cite{katz:rls} there is a smooth quasi-projective morphism of schemes 
\[ \pi: \mathrm{Hyp}\lra U \]
with geometrically connected fibres and an action of $\mu_N$ on $\mathrm{Hyp}$ over $U$, and a faithful character $\chi: \mu_N(\barQQ) \ra \barQQ_l^\times$, so that 
\[ \cF \cong (\mathrm{Gr}_r^W R^r\pi_! \barQQ_l)_{\barQQ}^\chi. \]
Here $\mathrm{Gr}^W$ denotes the weight filtration; $r$ denotes the relative dimension of $\mathrm{Hyp}$ over $U$; and the subscript $\barQQ$ indicates base change to $U_{\barQQ}$. More specifically, if $\infty \in S$, then $\mathrm{Hyp}$ is the hypersurface in $\G_m \times \Aff^{r+1}$ defined by the equation
\[ Y^N = \prod_{s \in S-\{\infty\}} \prod_{i=1}^{r+1} (X_i - s)^{e_i(s)} \prod_{i=1}^r (X_{i+1}-X_i)^{f_i}, \]
where $e_i(s)$ and $f_i \in \Z_{\geq 0}$ and none of the $f_i$'s is divisible by $N$. Moreover the action of $\mu_N$ is by multiplication on $Y$. Any choice of $r$, $e_i(s)$, $f_i$ and (faithful) $\chi$ can arise for some $\cF$.  

If $K$ is a number field and if $u \in U(K)$, then there is a $\mu_N$-equivariant projective compactification \footnote{Existence of such was announced but not written up by Hironaka. For a proof, see \cite{abramovich-wang}; namely, consider the product of all compositions with $g \in \mu_N(\barQQ)$ of a given $\mathrm{Hyp} \to \mathbb{P}^N$ to produce a $\mu_N$-equivariant inclusion of $Y$ as a quasi-projective sub-variety of some $\mathbb{P}^M$; take the projective closure and apply \cite[Theorem 0.1]{abramovich-wang} to the result.}
\[ \mathrm{Hyp}_u \into \overline{\mathrm{Hyp}}_u \]
with complement $D = \bigcup_i D_i$ a union of smooth divisors with normal crossings. Then by a standard argument (see e.g. \cite{deligne:hodgeI}) 
\[ (\mathrm{Gr}_r^W R^r\pi_! \barQQ_l)^\chi_u \cong \ker (H^r(\overline{\mathrm{Hyp}}_{u,\barK},\barQQ_l) \lra \bigoplus_{i \in I} H^r(D_{i,\barK},\barQQ_l) )^\chi. \]
If $K \supset \Q(\zeta_N)$, so that $\chi:\mu_N \ra \G_m$ over $K$, then we define a pure object $M(\cF,u)$ of $\cM_{K,\Q(\zeta_N)}$ to be the $\chi$ component of the  kernel of the map of motivated motives
\[ \mathrm{Gr}_r^W (\overline{\mathrm{Hyp}}_u,1,0) \lra \bigoplus_{i \in I} \mathrm{Gr}_r^W (D_i,1,0). \]
Then for any rational prime $l'$ and any embedding $i':\Q(\zeta_N) \into \barQQ_{l'}$ we have
\[ \begin{array}{rcl} \!\!\! H_{\lambda'}(M(\cF,u)) \otimes_{\Q(\chi)_{\lambda'},i'} \barQQ_{l'} & \!\!\!\! \cong & \!\!\!\! \mathrm{Gr}^W_r H^r_c(\mathrm{Hyp}_{u,\barK},\barQQ_{l'})^{i'(\chi)}\\
& \!\!\!\! \cong &\!\!\!\! \ker (H^r(\overline{\mathrm{Hyp}}_{u,\barK},\barQQ_{l'}) \lra \bigoplus_{i \in I} H^r(D_{i,\barK},\barQQ_{l'}) )^{i'(\chi)}, \end{array} \]
where $\lambda'$ denotes the prime of $\Q(\zeta_N)$ induced by $i'$. 
We do not claim that $M(\cF,u)$ only depends on $\cF$ and $u$. To the best of our knowledge it also depends on the choices of ${\mathrm{Hyp}}$ and $\chi$. It is however independent of the choice of compactification  $\overline{\mathrm{Hyp}}_u$.\footnote{Given two such equivariant projective compactifications we can find a third such compactification mapping to both of them. To see this one just needs to apply \cite[Theorem 0.1]{abramovich-wang} to the closure of the diagonal embedding of $\mathrm{Hyp}$ into the product of the two given compactifications. Then the motivated motive resulting from the third compactification maps to the motivated motives arising from the first two compatifications. Moreover, because they induce isomorphisms in cohomology, these maps are isomorphisms.}  {\em When we make an assertion about $M(\cF,u)$ it should be read as applying whatever auxiliary choices are made.}

It follows from \cite{pink} that for all but finitely many places $v$ of $K$, all $j$ sufficiently large and all $\zeta \in \mu_N(K)$ the alternating sum 
\[ \sum_s (-1)^s \tr (\zeta \Frob_v^j)|_{H^s_c(\mathrm{Hyp}_{u,\barK},\barQQ_{l'})} \]  
lies in $\Q$ and is independent of $l'$. Thus 
\[ \sum_s (-1)^s \tr \Frob_v^j|_{H^s_c(\mathrm{Hyp}_{u,\barK},\barQQ_{l'})^{i'(\chi)}} \]  
lies in $\Q(\zeta_N)$ and is independent of $l'$ and $i'$. By theorem 8.4.1(2) of \cite{katz:rls} we see that this latter sum simply equals
\[ (-1)^r \tr \Frob_v^j|_{H^r_c(\mathrm{Hyp}_{u,\barK},\barQQ_{l'})^{i'(\chi)}}. \]
Thus 
\[ \tr \Frob_v^j|_{\mathrm{Gr}^W_r H^r_c(\mathrm{Hyp}_{u,\barK},\barQQ_{l'})^{i'(\chi)}} \]
is independent of $l'$ and $i'$. We deduce that $M(\cF,u)$ is compatible.

\begin{lem} If $\Q(u)$ is totally real then 
\[ {}^cM(\cF,u) \cong M(\cF,u)^\vee(-r). \]
Thus $M(\cF,u)$ is compatible, polarized and pure.
\end{lem}

\pfbegin ${}^cM(\cF,u)$ is the $\chi\circ c = \chi^{-1}$ component of the  kernel of the map of motivated motives
\[ \mathrm{Gr}_r^W (\overline{\mathrm{Hyp}}_u,1,0) \lra \bigoplus_{i \in I} \mathrm{Gr}_r^W (D_i,1,0). \]

$M(\cF,u)^\vee (-r)$ is isomorphic to the $\chi$ component of the cokernel of the map
\[ \mathrm{Gr}_r^W (\overline{\mathrm{Hyp}}_u,1,0)^\vee(-r) \lla \bigoplus_{i \in I} \mathrm{Gr}_r^W (D_i,1,0)^\vee(-r). \]
By definition of the dual of an object in $\cM_{K,\Q(\chi)}$ this map can be thought of as a map
\[ \mathrm{Gr}_r^W (\overline{\mathrm{Hyp}}_u,1,0) \lla \bigoplus_{i \in I} \mathrm{Gr}_r^W (D_i,1,-1). \] 
Moreover, if $\zeta \in \mu_N$, then the action $\zeta^\vee$ of $\zeta$ on $(\overline{\mathrm{Hyp}}_u,1,0)^\vee(-r)$ corresponds to the action of $\zeta^{-1}$ on $(\overline{\mathrm{Hyp}}_u,1,0)$. (The transpose of the graph of an automorphism is the graph of the inverse austomorphism.) Thus $M(\cF,u)^\vee (-r)$ is isomorphic to the $\chi^{-1}$ component of the cokernel of the map
\[ \mathrm{Gr}_r^W (\overline{\mathrm{Hyp}}_u,1,0) \lla \bigoplus_{i \in I} \mathrm{Gr}_r^W (D_i,1,-1). \] 

There is a natural map from the kernel of
\[ \mathrm{Gr}_r^W (\overline{\mathrm{Hyp}}_u,1,0) \lra \bigoplus_{i \in I} \mathrm{Gr}_r^W (D_i,1,0) \]
to the cokernel of the map
\[ \mathrm{Gr}_r^W (\overline{\mathrm{Hyp}}_u,1,0) \lla \bigoplus_{i \in I} \mathrm{Gr}_r^W (D_i,1,-1) \] 
induced by the the identity on $(\overline{\mathrm{Hyp}}_u,1,0)$. 
To prove the lemma it suffices to show that after taking $\chi$-components this map is an isomorphism, or even that it is an isomorphism after applying $H_l$. Applying $H_l$ we have a commutative diagram with exact rows:
\[ \begin{array}{ccccccc}
\!\!(0) &\!\!  \ra &\!\! \gr^W_rH^r_c(\mathrm{Hyp}_{u,\barQQ},\Q_l) & \!\!\ra & \!\!H^r_c(\overline{\mathrm{Hyp}}_{u,\barQQ},\Q_l) & \!\!\ra & \!\!\bigoplus_{i \in I} H^r_c(D_{i,\barQQ},\Q_l) \\ && \da && || && \\
\!\!(0) &  \!\!\la & \!\!\gr^W_rH^r(\mathrm{Hyp}_{u,\barQQ},\Q_l) &\!\! \la & \!\!H^r(\overline{\mathrm{Hyp}}_{u,\barQQ},\Q_l) &\!\! \la & \!\!\bigoplus_{i \in I} H^{r-2}(D_{i,\barQQ},\Q_l(-1)), \end{array} \]
where the second row is the dual of the first row and where the first vertical map is the natural map. Thus it suffices to show that
\[  \gr^W_rH^r_c(\mathrm{Hyp}_{u,\barQQ},\Q_l)^\chi \lra  H^r(\mathrm{Hyp}_{u,\barQQ},\Q_l) \]
is injective.

To do this we return to the constructions of chapter 8 of \cite{katz:rls}. Looking at the proof of theorem 8.4.1 of \cite{katz:rls} we see that it suffices to check that, in the notation of theorem 8.3.5 of \cite{katz:rls}, the map 
\[ R^r(\mathrm{pr}_{r+1})_! \cL \lra R^r(\mathrm{pr}_{r+1})_* \cL \]
has image $\cH_r$. 

To prove this we first note that one can add to the conclusion of lemma 8.3.2 of \cite{katz:rls} the assertion that the natural map 
\[ \mathrm{NC}_\chi(\cF) \lra R^1 (\mathrm{pr}_2)_* (\mathrm{pr}_1^* (\cF) \otimes \cL_{\chi(X_2-X_1)}) \]
has image $\mathrm{MC}_\chi(\cF)$. (In the notation of that lemma.) Indeed this follows from the other conclusions of that lemma and the fact that
\[ R^1 (\overline{\mathrm{pr}}_2)_*(j_* (\mathrm{pr}_1^* (\cF) \otimes \cL_{\chi(X_2-X_1)})) \lra R^1 (\mathrm{pr}_2)_* (\mathrm{pr}_1^* (\cF) \otimes \cL_{\chi(X_2-X_1)})\]
is injective. (We remark that this statement is essentially the fact that for $j \colon V \subset Y$ an inclusion of an open $V$ in a smooth proper curve $Y$, and $\mc{F}$ a lisse sheaf on $V$, the parabolic cohomology 
\[ H^1(Y, j_{*}\mc{F}) \cong \Im( H^1_c (V, \mc{F}) \to H^1(V, \mc{F})).) \]

Now we return to our desired strengthening of theorem 8.3.5 of \cite{katz:rls}. We argue by induction on $r$. The case $r=1$ is just our strengthening of lemma 8.3.2 of \cite{katz:rls}. In general by the inductive hypothesis we know that
\[ \cH_{r-1} \lra R^{r-1}(\mathrm{pr}_{r})_* \cL \]
is injective, where we use the notation of the proof of theorem 8.3.5 of \cite{katz:rls}.
Pulling back to $\A(n,2)$ and again applying our strengthening of lemma 8.3.2 of \cite{katz:rls} we see that
\[ \cH_r \into R^1(\mathrm{pr}_2)_* R^{r-1}(\mathrm{pr}'_{r+1})_* \cL, \]
where $\mathrm{pr}'_{r+1}$ denotes the map $\A(n,r+1) \ra \A(n,2)$, which forgets the first $r-1$ $X$-coordinates. As $\mathrm{pr}_{r+1}=\mathrm{pr}_2 \circ \mathrm{pr}'_{r+1}$, we have a spectral sequence with second page
\[ E_2^{i,j}=R^i(\mathrm{pr}_2)_* R^{j}(\mathrm{pr}'_{r+1})_* \cL \Rightarrow R^{i+j} (\mathrm{pr}_{r+1})_* \cL. \]
However $\mathrm{pr}_2$ is affine of relative dimension $1$. Thus $R^i(\mathrm{pr}_2)_*$ vanishes for $i>1$ and our spectral sequence degenerates at $E_2$. We deduce that 
\[ \cH_r \into R^1(\mathrm{pr}_2)_* R^{r-1}(\mathrm{pr}'_{r+1})_* \cL \into R^{r} (\mathrm{pr}_{r+1})_* \cL, \]
as desired.
\pfend

If we want to apply corollary \ref{motpm} to $M(\cF,u)$ we must calculate
\[ \dim_{\overline{\Q(\chi)}} \gr^j H_\DR(M(\cF,u)) \otimes_{K \otimes \Q(\chi), \tau \otimes 1} \overline{\Q(\chi)} \]
for all $j$ and all $\tau: K \into \overline{\Q(\chi)}$. Choosing an embedding $\imath: \overline{\Q(\chi)} \into \C$ we see that this equals the dimension of the $\imath(\chi)$-component of $j^{th}$-graded piece (for the Hodge filtration) of the kernel of
\[ H^r(\overline{\mathrm{Hyp}}_u(\C), \C) \lra \bigoplus_{i \in I} H^r(D_i(\C),\C), \]
where the $\C$ points of these schemes over $K$ are calculated via $\imath \circ \tau: K \into \C$. This is the same as 
\[ \dim_\C \gr^j_F \gr_r^W H^r_c({\mathrm{Hyp}}_u(\C), \C)^{\imath(\chi)}, \]
where $\gr^j_F$ denotes the graded pieces for the Hodge filtration. This in turn equals
the rank of 
\[ \gr^j_F \gr_r^W (R^r\pi_! \C)^{\imath(\chi)}. \]
This does not depend on the choice of $u$. Extend $\imath$ to an isomorphism $\barQQ_l \iso \C$ and let $\imath'$ be an extension of $\imath \circ \tau$ to an embedding $\barQQ \into \C$. Then
\[  \gr_r^W (R^r\pi_! \C)^{\imath(\chi)} \cong \imath'_* (\cF \otimes_{\barQQ_l,\imath} \C)^\ann. \]
We see that for any $\imath: \barQQ_l \iso \C$ and $\imath': \barQQ \into \C$ the sheaf
\[ \cF_{\imath,\imath'}=\imath'_* (\cF \otimes_{\barQQ_l,\imath} \C)^\ann\]
 admits a polarizable variation of $\C$-Hodge structures. (For variations of Hodge structures we will follow the terminology of section 3.4 of \cite{sign}.)  According to proposition 1.13 of \cite{del:uni} this polarizable variation of $\C$-Hodge structures is unique up to translating the numbering. 
We will say that $\cF$ is {\em regular} if for all $\imath, \imath'$ and all $j$ we have
\[ \dim_\C \gr^j_F \cF_{\imath,\imath'} \leq 1. \]
Then we deduce the following consequence of corollary \ref{motpm}.

\begin{cor}\label{rmpm} Suppose that $S$ is a finite subset of $\PP^1(\Q)$ and that $\cF$ is a regular rigid local system on $(\PP^1-S)_{\barQQ}$ with quasi-unipotent monodromy of order dividing $N$. Suppose that $F$ is a CM field containing a primitive $N^{th}$ root of unity and that $u \in F^+$. Then there is a finite Galois CM extension $F'/F$ and a decomposition (perhaps after extending the field of coefficients)
\[ \CH(M(\cF,u)) \equiv \CR_1 \oplus ... \oplus \CR_s \]
such that each $\CR_i|_{G_{F'}}$ is automorphic. 

In particular $\CH(M(\cF,u))$ is strictly pure and, if $\imath:\Q(\zeta_N) \into \C$, then the completed L-function $\Lambda(\imath \CH(M(\cF,u)),s)$ has meromorphic continuation to the whole complex plane and satisfies the functional equation
\[ \Lambda(\imath \CH(M(\cF,u)),s)=\epsilon(\imath \CH(M(\cF,u)),s) \Lambda(\imath \CH(M(\cF,u))^\vee, 1-s). \]
\end{cor}

To apply this corollary one must be able to calculate the
\[ \dim_\C \gr^j_F \cF_{\imath,\imath'}. \]
This is discussed in \cite{dettweiler-sabbah}, but seems in general to be a complicated question.  
We will discuss a more explicit condition that implies regularity, but which is probably a much stronger condition. 

We will say that $\cF$ {\em has somewhere maximally quasi-unipotent monodromy} if for some $s\in S$ the monodromy of $\cF$ at $s$ has only one Jordan block. We have the following lemma.

\begin{lem} If $S$ is a finite subset of $\PP^1(\Q)$ and if $\cF$ is a rigid local system on $(\PP^1-S)_{\barQQ}$ with quasi-unipotent monodromy of order dividing $N$ and with somewhere maximally quasi-unipotent monodromy, then $\cF$ is regular. \end{lem}

\pfbegin 
Choose $s \in S$ such that the monodromy of $\cF$ at $s$ has only one Jordan block. 
Twisting by the pull back of the lisse sheaf on $\G_{m,\barQQ}$ corresponding to a finite order character of $\pi_1(\G_{m,\barQQ})$ under some linear map $(\PP^1-S) \ra \G_m$,  we may suppose that $\cF$ has unipotent monodromy at $s$.

There is a limit mixed Hodge structure $\cF_{\imath,\imath',s}$ as follows: $\cF_{\imath,\imath',s}$ is a $\C$-vector space with
\begin{itemize}
\item two decreasing exhaustive and separated filtrations $\Fil^i_F$ and $\barFil^i_{F}$,
\item an increasing exhaustive and separated filtration $\Fil^W_j$,
\item and a nilpotent endomorphism $N$
\end{itemize}
with the following properties:
\begin{enumerate}
\item $\Fil^W_j=\sum_{i_1+r=1+j+i_2} (\ker N^{i_1}) \cap (\Im N^{i_2})$;
\item $\cF_{\imath,\imath',s}$ with its automorphism $\exp N$ is isomorphic to $\cF_{\imath,\imath',z}$ with its monodromy operator for any $z \in \PP^1(\C)-S$;
\item $\dim_\C \gr_F^i \cF_{\imath,\imath',s}= \dim_\C  \gr_F^i\cF_{\imath,\imath',z}$ and $\dim_\C \bargr_F^i \cF_{\imath,\imath',s}= \dim_\C  \bargr_F^i\cF_{\imath,\imath',z}$ for any $z \in \PP^1(\C)-S$;
\item $\Fil^i_F$ and $\barFil_F^i$ induce on $\mathrm{gr}^W_j \cF_{\imath,\imath',s}$ a pure $\C$-Hodge structure of weight $j$;
\item\label{trans} we have equalities $N\Fil^i_F \cF_{\imath,\imath',s} = (\Im N) \cap \Fil^{i-1}_F \cF_{\imath,\imath',s}$ and
$N \barFil^i_F \cF_{\imath,\imath',s}= (\Im N) \cap \barFil^{i-1}_F \cF_{\imath,\imath',s}$.
\end{enumerate}
This is true for any polarizable variation of pure $\C$-Hodge structures with unipotent monodromy on an open  complex 
disc minus one point. This follows from the corresponding fact for polarizable variations of pure $\R$-Hodge structures, by the dictionary between variations of $\C$-Hodge structures and variations of $\R$-Hodge structures with an action of $\C$. (See for instance section 3.4 of \cite{sign}.) In the case of variations of $\R$-Hodge structures it follows from theorem 6.16 of \cite{schmid}. Note that theorem 6.16 of \cite{schmid} seems to make the assumption that the variation of Hodge structures is the cohomology of a family of smooth projective varieties. However as pointed out in \cite{cos}, see in particular the first paragraph of page 462 of that paper, this plays no role in the proof.

From assertion (\ref{trans}) above we see that for all $i$ we have
\begin{equation}\label{inj} N: \Fil^i_F \cF_{\imath,\imath',s} /(\Fil^{i+1}_F\cF_{\imath,\imath',s}+\ker N) \into \gr_F^{i-1}\cF_{\imath,\imath',s}. \end{equation}
Choose $i_0$ maximal such that $\ker N \subset \Fil^i_F\cF_{\imath,\imath',s}$. For $i< i_0$ we have
\[ N : \gr^i_F\cF_{\imath,\imath',s} \into \gr^{i-1}_F\cF_{\imath,\imath',s} \]
and so in fact $\gr^i_F\cF_{\imath,\imath',s}=(0)$. Thus $\Fil^{i_0}_F \cF_{\imath,\imath',s}=\cF_{\imath,\imath',s}$. 
As $\ker N$ is one dimensional we see that 
\[ \Fil^{i_0+1}_F \cap \ker N = (0). \]
Because of the injection \ref{inj} we see that
\[ \Fil^{i_0}_F\cF_{\imath,\imath',s}=\Fil^{i_0+1}_F \cF_{\imath,\imath',s}\oplus \ker N, \]
and that, for $i\geq i_0$,
\[ N^{i-i_0}: \gr^i_F \cF_{\imath,\imath',s}\into \gr^{i_0}_F\cF_{\imath,\imath',s}. \]
Thus for all $i$
\[ \dim_\C \gr_F^i\cF_{\imath,\imath',s} \leq 1, \]
and the lemma follows.
\pfend 

One can use the constructions of section 8.3 of \cite{katz:rls} to construct many examples of rigid local systems $\cF$ on $U_{\barQQ}$ with quasi-unipotent monodromy of order dividing $N$ and with somewhere maximally quasi-unipotent monodromy. More classically, Corollary \ref{rmpm} applies to any hypergeometric local system on $\mathbb{P}^1-\{0, 1, \infty\}$ with somewhere maximally unipotent local monodromy. By \cite[Theorem 3.5]{beukers-heckman}, we get such an example for any choice of function $m: \mu_\infty(\barQQ)-\{1\} \ra \Z_{\geq 0}$ with finite support. It will have rank $\sum_\zeta m(\zeta)$. The local monodromies at the three punctures will be a quasi-reflection, a single unipotent Jordan block, and a Jordan form where each $\zeta \in \mu_\infty(\barQQ)-\{1\}$ appears in a single Jordan block with length equal to $m(\zeta)$. Finally, here is a less classical, non-hypergeometric example, generalizing a construction in \cite{dettweiler-reiter:G2rls}: 

If $\chi$ is a continuous $\barQQ_l$-valued character of $\pi_1(\G_{m,\barQQ})$ let $\cL_\chi$ denote the corresponding lisse sheaf on $\G_{m,\barQQ}$. If $\chi_1$ and $\chi_2$ are continuous $\barQQ_l$-valued characters of $\pi_1(\G_{m,\barQQ})$ let $\cL(\chi_1,\chi_2)$ denote the lisse sheaf on $(\PP^1-\{0,1,\infty\})_{\barQQ}$ which is the tensor product of the pull back of $\cL_{\chi_1}$ under the identity map with the pull-back of $\cL_{\chi_2}$ under the map $t \mapsto t-1$. We will write $\mathbf{1}$ for the trivial character of $\pi_1(\G_{m,\barQQ})$ and $-\mathbf{1}$ for the unique character of $\pi_1(\G_{m,\barQQ})$ of exact order $2$. Define lisse sheaves $\cF_i$ on $(\PP^1-\{0,1,\infty\})_{\barQQ}$ recursively by
\[ \cF_0 = \cL(-\mathbf{1},-\mathbf{1}) \]
and
\[ \cF_{2i-1}= \cL(\mathbf{1},-\mathbf{1}) \otimes {\rm MC}_{-\mathbf{1}}(\cF_{2i-2}) \]
and
\[ \cF_{2i}= \cL(-\mathbf{1},\mathbf{1}) \otimes {\rm MC}_{-\mathbf{1}}(\cF_{2i-1}), \]
for $i \in \Z_{>0}$. Here ${\rm MC}_{-\mathbf{1}}$ denotes the functor associated to the representation $-\mathbf{1}$ of $\pi_1(\G_{m,\barQQ})$ as described in section 8.3 of \cite{katz:rls}. (It is closely related to the `middle convolution'.)  It follows inductively, using section 5.1 of \cite{katz:rls} and proposition 1.2.1 of \cite{dettweiler-reiter:G2rls}, that each $\cF_i$ is a rigid local systems $\cF$ on $(\PP^1-\{0,1,\infty\})_{\barQQ}$ with quasi-unipotent monodromy of order dividing $2$ and with somewhere maximally quasi-unipotent monodromy, and that its full monodromy is give by the following table. (Knowing the monodromy of $\cF_i$ everywhere is needed to calculate the monodromy of $\cF_{i+1}$ at $\infty$. This example is a generalization of the special case $\cF_6$ which is considered in \cite{dettweiler-reiter:G2rls}, and which is a $G_2$-local system.) 
\begin{table}[!h]
\centering

\begin{tabular}{c|c|c|c}
 & at $0$ & at $1$ & at $\infty$ \\ \hline
 $i \equiv 0 \mod 4$ & $\mathbf{1}^{\oplus \frac{i}{2}} \oplus (-\mathbf{1})^{\oplus(\frac{i}{2}+1)}$ & $-\mathbf{1} \oplus U(2)^{\oplus \frac{i}{2}}$ & $U(i+1)$ \\  \hline
$i \equiv 1 \mod 4$ & $U(2)^{\oplus \frac{i+1}{2}}$ & $(-\mathbf{1} \otimes U(2)) \oplus (-\mathbf{1})^{\oplus \frac{i-1}{2}} \oplus \mathbf{1}^{\oplus \frac{i-1}{2}}$ & $U(i+1)$ \\ \hline
  $i \equiv 2 \mod 4$ & $\mathbf{1}^{\oplus \frac{i}{2}} \oplus (-\mathbf{1})^{\oplus(\frac{i}{2}+1)}$ & $U(3) \oplus U(2)^{\oplus \frac{i-2}{2}}$ & $U(i+1)$ \\ \hline
  $i \equiv 3 \mod 4$ & $U(2)^{\oplus \frac{i+1}{2}}$ & $U(2) \oplus \mathbf{1}^{\oplus \frac{i-3}{2}} \oplus (-\mathbf{1})^{\oplus \frac{i+1}{2}}$ & $U(i+1)$ \\ \hline

\end{tabular}
\caption{The local monodromies for the sheaves $\cF_i$}
\label{table}
\end{table}
In particular corollary \ref{rmpm} applies to each of the sheaves $\cF_i$. (Note that because of the regular quasi-unipotent monodromy at $\infty$, this result would previously have been available for $u \in F-\cO_F$, because in this case $H_\lambda(M(\cF_i,u))$ is automatically irreducible. See \cite{blggt}.)

\end{document}